\documentclass{article}
\usepackage{amssymb}
\usepackage{amsthm}
\usepackage{amsmath}
\usepackage[cp1250]{inputenc}

\newtheorem{theorem}{Theorem}[section]
\newtheorem{lemma}[theorem]{Lemma}

\newtheorem{corr}[theorem]{Corrolary}
\newtheorem{definition}[theorem]{Definition}
\theoremstyle{definition}

\newcommand{\comment}[1]{}

\def\qed{{\hfill{\vrule height5pt width3pt depth0pt}\medskip}}

\begin{document}
\begin{center}
{\Large \bf Hyperbolicity and  Averaging for the Srzednicki-W{\'o}jcik equation}

 \vskip 0.5cm
Piotr~Zgliczy\'nski\footnote{Research has been supported by Polish
National Science Centre grant 2015/19/B/ST1/01454}

 \vskip 0.5cm

{\small Institute of Computer Science,
Faculty of Mathematics and Computer Science, Jagiellonian University}\\
{\small  ul. S. {\L}ojasiewicza 6, 30-348 Krak\'ow, Poland}

\vskip 0.5cm

umzglicz@cyf-kr.edu.pl, piotr.zgliczynski@ii.uj.edu.pl

\vskip 0.5cm

\today
\end{center}

\begin{abstract}
  For the Srzednicki-W{\'o}jcik equation, the planar nonautonomous ODE parameterized by $\kappa \in \mathbb{R}$,
  $$
      z'=\overline{z}(1+ |z|^2 \exp(i \kappa t)), \qquad z(t) \in \mathbb{C}
  $$
  using averaging  we show how the region
  of hyperbolicty grows with $|\kappa|$. Based on this we give bounds on the sizes of bounded orbits.
\end{abstract}


\section{Introduction}

%The main topic of this paper is the averaging of fast oscillations in ODEs on relatively large domains.
In this paper we investigate some aspect of dynamics of Srzednicki-W{\o}jcik equation  \cite{SW}
\begin{equation}
  z'=\overline{z}(1+ |z|^2 \exp(i \kappa t)), \qquad z(t) \in \mathbb{C}. \label{eq:srzednicki}
\end{equation}
This equation has been  studied by topological methods \cite{SW,WZ,OW,MS}, with the goal of establishing
 the existence of symbolic dynamics. Main tool was the method of isolating segments introduced in \cite{S1,SW} with papers \cite{WZ,OW} improving estimates
 so finally the existence of symbolic dynamics was established for $0<\kappa <0.5044$. In \cite{MS} using related ideas a computer assisted proof was
 obtained for $0< \kappa \leq 1$.

For equation (\ref{eq:srzednicki})  we establish the following results
\begin{itemize}
\item[1.] all solutions of (\ref{eq:srzednicki}) with sufficiently large $\|z(0)\| (\geq O(|\kappa|^{1/2}))$ go to infinity in finite time either forward or backward in time. This is the content of Theorem~\ref{lem:toInfinityInFiniteTime} in Section~\ref{sec:srzed-eq}.
\item[2.] in the neighborhood of $0$ there is a hyperbolic region for the dynamics induced by (\ref{eq:srzednicki}), which  contains a ball of size $O(|\kappa|^{1/4})$ i.e. it grows with $|\kappa|$. This is a consequence  of Theorem~\ref{thm:eq-srzed-hyperbolic} in Section~\ref{sec:srzed-eq}.
\item[3.] we give  a lower bound $O(|\kappa|^{1/4})$ and an upper bound $O(|\kappa|^{1/2})$ of the norm  for any bounded nonzero orbits. This is discussed in Section~\ref{subsec:srzeeq-bnds}.
\end{itemize}
Items 2 and 3 depend on developing explicit estimates for the difference between solutions of
\begin{equation}
z'=\overline{z}  \label{eq:auto-srzed}
\end{equation}
and (\ref{eq:srzednicki}), and  also for the solutions of corresponding variational equations on large domains.

We treat (\ref{eq:srzednicki}) as the perturbation of   problem (\ref{eq:auto-srzed}).
Since (\ref{eq:srzednicki}) is not a small perturbation of (\ref{eq:auto-srzed}) we cannot apply to the vector
fields (\ref{eq:srzednicki}) or (\ref{eq:auto-srzed})
topological tools such as the isolating blocks or isolating segments originating from
the Conley index theory \cite{MM,SW}, or the continuation based on the implicit function
theorem. Instead we use the averaging method to obtain bounds on the difference of solutions
of these systems.

  Our approach to averaging differs from the dominant one \cite{AKN97,Nei84,GH,Hale,SV,JS,L}, which following the pioneering work of (\cite{BM,BZ}) used  suitable coordinates
changes to control the influence of the rapidly oscillating perturbation.

In this work we directly estimate the influence of the rapidly oscillating perturbation
on the solution, as it was proposed in the PDE context by Henry \cite{He} (who attribute this approach in the ODE context to Gikhman).
  We exploit the phenomenon  known in the numerical analysis community
\cite{I}, that fast oscillations are your friend, if one integrates them
first. In this work following this advice we consider shifts along
the trajectory by a small time step. It turns out that the time shift
map for (\ref{eq:srzednicki}) is a small perturbation of the time shift
of (\ref{eq:auto-srzed}) if  $|\kappa| \to \infty$.
In literature on averaging for ODEs somewhat analogous methods can be found in \cite{BL,Wa}. In \cite{Wa} an idea is developed that in a suitable weak topology
the averaged equation is the limit of the problem with fast oscillation with the oscillation frequency going to infinity. In a sense this observation
underlines our approach.

In contrast to other works on averaging we insist on explicit estimates. In the context of (\ref{eq:srzednicki}) this requires obtaining a priori bounds
on solutions, to obtain  the sets on which we can compute difference between solutions of (\ref{eq:srzednicki}) and (\ref{eq:auto-srzed}) based on the averaging.

The approach has been originally developed with the aim to be applicable to dissipative PDEs. We applied it to viscous Burgers equation and
Navier-Stokes equations with periodic boundary conditions
in paper \cite{CyZ} and in \cite{CMTZ} to Navier-Stokes equations and the damped-Euler equation.

\subsection{Notation}
Consider   nonautonomous ODE
\begin{equation}
  x'(t)=f(t,x(t)), \label{eq:notation-Non_auto}
\end{equation}
 where $x \in \mathbb{R}^n$ and $f: \mathbb{R} \times \mathbb{R}^n \to \mathbb{R}^n$ is regular enough to guarantee
 the
uniqueness of the initial value problem $x(t_0)=x_0$.  We set
 $\varphi(t_0,t,x_0)=x(t_0+t)$, where $x(t)$ is a solution of (\ref{eq:notation-Non_auto}) with  initial condition  $x(t_0)=x_0$.
 Obviously in each context it will be clearly stated what is the ordinary differential equation generating $\varphi$. We will
 sometimes refer to $\varphi$ as to the local process generated by  (\ref{eq:notation-Non_auto}).

 For  matrix $U$ by $U^t$ we will denote its transpose. For a square matrix $U$ we will denote its spectrum by $\mbox{Sp}(U)=\{ \lambda \in \mathbb{C} \ | \ \lambda \ \mbox{is an eigenvalue of $U$}\}$. If $A \in \mathbb{R}^{n\times n}$ by $\mu(A)$ we will denote its logarithmic norm, which is defined
 by
 \begin{equation}
   \mu(A) = lim_{h \to 0^+} \frac{\|I + hA\|- 1}{h}.
 \end{equation}
 For the properties of logarithmic norm and its relation with the Lipschitz constant
 for the flow induced by ODEs see \cite{HNW,KZ} and literature cited there.
 The logarithmic norm depends on the norm used. We will always assume that we are
 using the euclidean norm.

%Quite often, in connection with bounds based on the logarithmic norms,   the following expression will appear
%$h_l(t)=\frac{e^{lt}-1}{l}$, where $l$ is a fixed parameter. For
%$l=0$ we set $h_0(t):=\lim_{l \to 0}\frac{e^{lt}-1}{l}=t$

  For a function of several variables we will often use $D_tf$ and $D_z f$
to denote the partial derivatives. For example, $D_z
f(t,z)=\frac{\partial f}{\partial
 z}(t,z)$.

% to jest to samo co w rapid-osc-dla-pdes i tak chce to utrzymac

\section{Averaging - basic estimates} % 
\label{sec:basic-estm}

The results from Sections~\ref{subsec:lin-non-auto} and~\ref{subsec:non-prob1} are also contained in \cite{CyZ}. We include them here for the sake of making the paper reasonably self-contained.

\subsection{Linear nonautonomous equations}
\label{subsec:lin-non-auto}

Assume that $A: \mathbb{R} \to \mathbb{R}^{n \times n}$ is continuous and
for $k=1,\dots,m$, $v_k:\mathbb{R}  \to \mathbb{R}$
are $C^1$, $g_k:\mathbb{R} \to \mathbb{R}$ are continuous.

Let us consider the following non-autonomous non-homogenous linear ODE
\begin{equation}
  x'(t)=A(t)x(t) + \sum_{k \in J_V}  g_k(\omega_k t) v_k(t),  \quad x \in \mathbb{R}^n.  \label{eq:linNonHomoNonAuto}
\end{equation}
The set $J_V$ in the sum in (\ref{eq:linNonHomoNonAuto}) might be  finite or infinite, or the sum might be an integral over some measure on $J_V$.

For each $k \in J_V$ let $G_k(t)$ be a primitive of $g_k$, so
\begin{equation}
 G_k'(t)=g_k(t).
\end{equation}
We will assume later that $|G_k(t)|$'s are  bounded. This is the reflection of the oscillating nature of $g_k$.

Let $M(t,t_0)$ be a fundamental matrix of solutions of the homogenous version of (\ref{eq:linNonHomoNonAuto})
\begin{equation}
   x'(t)=A(t)x(t). \label{eq:linNonAuto}
\end{equation}
This means that for any $t_0 \in \mathbb{R}$ and $x_0 \in \mathbb{R}^n$ the function $x(t)=M(t,t_0)x_0$ solves
(\ref{eq:linNonAuto}) with the initial condition $x(t_0)=x_0$.

It is well known that $M$ has the following properties
\begin{eqnarray}
 M(t_0,t_0)&=&I, \\
 M(t,t_0)^{-1}&=&M(t_0,t), \\
 \frac{\partial}{\partial t} M(t,t_0)&=&A(t) M(t,t_0), \\
  \frac{\partial}{\partial t_0} M(t,t_0)&=&- M(t,t_0) A(t_0). \label{eq:Mfund-4}
\end{eqnarray}

 The general solution of (\ref{eq:linNonHomoNonAuto}) is given by
\begin{eqnarray}
 \varphi(t_0,t,x_0)&=&M(t_0+t,t_0)x_0 +   \label{eq:int-sol-fund-matrix} \\
   & & \int_{0}^{t} M(t_0+t,t_0+s)  \sum_{k \in J_V} g_k(\omega_k (t_0+s)) v_k(t_0+s)
   ds. \nonumber
\end{eqnarray}

We compute the integral in the above formula as follows. Using the integration by parts and (\ref{eq:Mfund-4}) we obtain for $k\in J_V$
\begin{eqnarray}
 I_k(t+t_0):=\int_{0}^{t}  g_k(\omega_k(t_0+s)) M(t_0+t,t_0+s) v_k(t_0+s)ds=  \label{eq:Dk-lin} \\
 \left. \frac{G(\omega_k(t_0+s))}{\omega_k} M(t_0+t,t_0+s)v_k(t_0+s)\right|^{s=t}_{s=0} +  \nonumber \\
-\frac{1}{\omega_k} \int_{0}^t G_k(\omega_k(t_0+s)) \frac{\partial }{\partial s}\left(M(t_0+t,t_0+s) v_k(t_0+s)\right)ds= \nonumber \\
 \frac{1}{\omega_k}(G_k(\omega_k(t_0+t)) v_k(t_0+t)- G_k(\omega_k t_0)M(t_0+t,t_0)v_k(t_0)) + \nonumber \\
 \frac{1}{\omega_k} \int_{0}^t G_k(\omega_k(t_0+s)) M(t_0+t,t_0+s)A(t_0+s)v_k(t_0+s)ds + \nonumber \\
- \frac{1}{\omega_k} \int_{0}^t G_k(\omega_k(t_0+s))
M(t_0+t,t_0+s)v'_k(t_0+s)ds \nonumber
\end{eqnarray}

For Galerkin projections of dissipative PDEs, while $\|A\|$ will not have any uniform bound independent of the projection dimension, we expect $\|Av_k\|$ to be uniformly bounded.

Therefore, we have proved that for the process generated by (\ref{eq:linNonHomoNonAuto}) it holds that
\begin{equation}
   \varphi(t_0,t,x_0)=M(t_0+t,t_0)x_0 + \sum_{k \in J_V} I_k(t+t_0).
\end{equation}

\subsection{Estimates for nonlinear problem}
\label{subsec:non-prob1}

Assume that $F: \mathbb{R} \times \mathbb{R}^{n} \to \mathbb{R}^n$ is $C^1$ function and
for $k=1,\dots,m$, $v_k:\mathbb{R} \times \mathbb{R}^n \to \mathbb{R}$
are $C^1$, $g_k:\mathbb{R} \to \mathbb{R}$ are continuous.

Consider problem
\begin{equation}
  x' = \tilde{F}(t,x):= F(t,x) + \sum_{k \in J_V}  g_k(\omega_k t)  v_k(t,x)\label{eq:prob1-loc}
\end{equation}
and its oscillation-free version
\begin{equation}
  y' = F(t,y).  \label{eq:auto-prob1-loc}
\end{equation}

%Let $\varphi(t_0,t,x)$ be the local process induced by
%(\ref{eq:prob1-loc}) and  $\varphi_a(t_0,t,x)$ is the local
%process induced by (\ref{eq:auto-prob1-loc}).

\begin{lemma}
\label{lem:prob1-estm-norm}
Let  $x:[t_0,t_0 + h] \to \mathbb{R}^n$ and $y:[t_0,t_0+h] \to \mathbb{R}^n$ be solutions to (\ref{eq:prob1-loc}) and (\ref{eq:auto-prob1-loc}),
respectively, such that $x(t_0)=y(t_0)$.

Let $W$ be a  compact set, such that for any  $t \in [0,h]$, the segment joining $x(t_0+t)$ and $y(t_0+t)$ is contained in $W$.

Assume that for $k \in J_V$ $G_k'(t)=g_k(t)$.

Assume that there exist constants $l$, $C(\dots)$  such that for all $k \in J_V$ it holds that
\begin{eqnarray}
   \sup_{t \in \mathbb{R}}\|G_k(t)\| &=& C(G_k), \\
   \sup_{t \in \mathbb{R}}|g_k(t)| &=& C(g_k), \\
  \sup_{z \in W, s \in [t_0,t_0+h]} \mu(D_z F(s,z))&=&l  \\
  \sup_{z \in W, s \in [t_0,t_0+h]} \|v_k(s,z)\|&=&C(v_k)  \\
  \sup_{z,z_1 \in W, s \in [t_0,t_0+h]} \|(D_z F(s,z)) v_k(s,z_1) \| &=& C\left(D_z F v_k\right) \\
  \sup_{z \in W, s \in [t_0,t_0+h]} \left\|\frac{\partial v_k}{\partial t} (s,z) \right\| &=&  C\left(\frac{\partial v_k}{\partial t}\right) \\
 \sup_{z \in W, s \in [t_0,t_0+h]} \left\|(D_z v_k (s,z)) \tilde{F}(s,z)  \right\| &=&  C(D_z v_k  \tilde{F} )
\end{eqnarray}

Assume that
\begin{eqnarray}
  \sum_{k\in J_V} C(G_k) C(v_k) < \infty  \label{eq:rosc-conv1} \\
    \sum_{k\in J_V} C(G_k) C\left(D_z F v_k\right) < \infty \\
     \sum_{k\in J_V} C(G_k) C\left(\frac{\partial v_k}{\partial t}\right) < \infty \\
      \sum_{k\in J_V} C(G_k) C\left(D_z v_k \tilde{F}\right) < \infty. \label{eq:rosc-conv4}
\end{eqnarray}

Then for $t \in [0,h]$ it holds that
\begin{equation}
 \|x(t_0+t)-y(t_0+t)\| \leq \sum_{k\in J_V} \frac{1}{|\omega_k|} b_k(t) %\leq \frac{1}{\inf_{k \in J_V} |\omega_k|}\sum_{k\in J_V} b_k(t) < \infty
   \label{eq:estm-diff-osc}
 \end{equation}
where  continuous functions $b_k:[0,h] \to \mathbb{R}_+$ depend on constants  $l$, $C(g_i)$, $C(G_i)$, $C(v_i)$, $C\left(D_z F v_i\right)$, $C\left(\frac{\partial v_i}{\partial t}\right)$ and
 $C\left(D_z v_i\tilde{F}\right)$ as follows
 \begin{eqnarray}
    b_k(t)&=&C(v_k) C(G_k)(1+e^{lt}) + C\left(D_z F v_k\right) C(G_k)(e^{lt}-1)/l +  \label{eq:bk-expression} \\
    & & C(G_k) \left(C\left(\frac{\partial v_k}{\partial t}\right) + C\left(D_z v_k \tilde{F}\right) \right)
(e^{lt}-1)/l \nonumber
 \end{eqnarray}

%\textbf{PZ: zalozenia (\ref{eq:rosc-conv1}--\ref{eq:rosc-conv4}) sa uzyte tylko do ponizszego stwierdzenia. Nie wydaje sie aby byly one warunkiem koniecznym,
%zbieznosc moglaby byc gwarantowana przez czynniki $\frac{1}{\omega_k}$ - popatrzec pozniej czy cos z tym zrobic. }

If $\inf_{k \in J_V} |\omega_k| >0$, the sum in (\ref{eq:estm-diff-osc}) is convergent.

\end{lemma}

\noindent
\textbf{Proof:}
Let $z(t)=x(t)-y(t)$. We have
\begin{eqnarray*}
 z'(t)=F(x(t)) - F(y(t)) +  \sum_{k\in J_V}  g_k(\omega_k t) v_k(t,x) =   \label{eq:prob1-diff}\\
   \left(\int_0^1 D_x F(t,s(x(t)-y(t)) + y(t))ds\right) \cdot z(t) + \sum_{k\in J_V}  g_k(\omega_k t) v_k(t,x(t)).
\end{eqnarray*}
Therefore
\begin{equation}
  z'(t) =  A(t)z(t) + \sum_{k\in J_V}  g_k(\omega_k t) v_k(t,x(t)),
\end{equation}
where
\begin{equation}
  A(t)=  \left(\int_0^1 D_x F(t,s(x(t)-y(t)) + y(t))ds\right).
\end{equation}
Let $M(t_1,t_0)$ is the fundamental matrix of solutions for $x'=A(t)x$.

Since $z(t_0)=0$, then from (\ref{eq:int-sol-fund-matrix}) and (\ref{eq:Dk-lin}) it follows that
\begin{equation}
  z(t_0+t)= \sum_{k\in J_V} I_k(t+t_0)
\end{equation}
where  %\textbf{PZ: czemu moglem zmienic kolejnosc sumowania?}
\begin{eqnarray*}
  I_k(t+t_0)=\frac{1}{\omega_k}(G_k(\omega_k(t_0+t)) v_k(t_0+t,x(t+t_0))- G_k(\omega_k t_0) M(t_0+t,t_0)v_k(t_0,x(t_0)) + \\
 \frac{1}{\omega_k} \int_{0}^t G_k(\omega_k(t_0+s)) M(t_0+t,t_0+s)A(t_0+s)v_k(t_0+s,x(t_0+s))ds + \\
- \frac{1}{\omega_k} \int_{0}^t G_k(\omega_k(t_0+s)) M(t_0+t,t_0+s)\left(\frac{d}{ds}v_k(t_0+s,x(t_0+s))\right)ds
\end{eqnarray*}

From the standard estimate for the
logarithmic norms (see for example Lemma 4.1 in \cite{KZ}) we know that
for $t \geq 0$ it holds that
\begin{equation*}
  \|M(t+t_0,t_0)\| \leq \exp(l t).
\end{equation*}
Hence, we obtain the following estimate of $I_k(t)$, for $t \in [0,h]$ and $k\in J_V$
\begin{eqnarray*}
  |\omega_k| \cdot \|I_k(t+t_0)\| \leq C(v_k) C(G_k)(1+e^{lt}) + C(D_z F v_k ) C(G_k)  \int_0^t e^{l(t-s)}ds + \\
     C(G_k) \left(C\left(\frac{\partial v_k}{\partial t}\right) + \sup_{s \in [0,h]} \left\|\frac{\partial v_k}{\partial z}(t_0+s,x(t_0+s)) x'(t_0+s)\right\|  \right) \int_0^t e^{l(t-s)}ds \leq \\
  C(v_k) C(G_k)(1+e^{lt}) + C(D_z F v_k) C(G_k)(e^{lt}-1)/l + \\
    C(G_k) \left(C\left(\frac{\partial v_k}{\partial t}\right) + C\left(D_z v_k \tilde{F}\right) \right)
(e^{lt}-1)/l
\end{eqnarray*}

This proves (\ref{eq:bk-expression}).

To finish the proof observe that if $|\omega_k| > \epsilon>0$ for all $k \in J_V$, then assumptions (\ref{eq:rosc-conv1}--\ref{eq:rosc-conv4}) together with the formula (\ref{eq:bk-expression}) imply the convergence of the sum in (\ref{eq:estm-diff-osc}).
\qed

\subsection{Estimates for higher order variational equations and persistence of normally hyperbolic objects}

The goal of this subsection it to state without much of precision what are the dynamical consequences of the just proved estimates.

The differentiation of (\ref{eq:prob1-loc}) with respect to initial conditions gives ODEs describing the evolution of the partial derivatives
of the local process $\varphi$. These equations have the same structure as the original ODEs, i.e.  consist from the part originating from
(\ref{eq:auto-prob1-loc}) and the rapidly oscillating one.
It is clear that under suitable conditions on higher derivatives of $v_k$ added to the assumptions of Lemma~\ref{lem:prob1-estm-norm} we can obtain the bounds of the following type
\begin{equation}
  \left\| \frac{\partial^{|\alpha|} \varphi}{\partial x^{\alpha}} (t_0,t,x) - \frac{\partial^{|\alpha|} \varphi_a}{\partial x^{\alpha}} (t_0,t,x)   \right\|
     \leq \frac{B_{|\alpha|}(t)}{\omega}
\end{equation}
where   $\alpha=(\alpha_1,\dots,\alpha_n) \in \mathbb{N}$, $|\alpha|=\alpha_1 + \dots + \alpha_n$, $\frac{\partial^{|\alpha|} \varphi}{\partial x^{\alpha}}$
is the partial derivative  $\frac{\partial^{|\alpha|} \varphi}{\partial x_1^{\alpha_1} \partial x_2^{\alpha_2}\cdots \partial x_n^{\alpha_n}}$ and
 $B_{|\alpha|}:[0,h] \to \mathbb{R}$ is continuous   function.

In particular from these conditions it will follow that the time shift map by $h$ along the trajectories  for $\varphi$ and $\varphi_a$ are $C^k$-close.

Assume that there exists, $M$, a compact normally hyperbolic invariant manifold (NHIM) for $\varphi_a$ (see \cite{HPS}). When seen in the extended phase space it gives rise to $\tilde{M}$
another NHIM, which is non compact but various important estimates are uniform on $\tilde{M}$. For example, a hyperbolic periodic orbit
for $\varphi_a$, which is topologically a circle, in the extended phase space becomes an infinite cylinder. Now we pass from the ODEs to the time shifts viewed in the extended phase space. As we said before these maps are $C^k$-close and one of them has NHIM  $\tilde{M}$.
 Despite the lack of compactness, due to uniform estimates (see for exmaple methods from \cite{CaZ}), one can show that HHIM $\tilde{M}$ exists also for the time shift
by $h$ for the process $\varphi$. Moreover, if the forcing is time periodic, then we obtain a true NHIM in the extended phase
space, which is periodic in time.

These observations are essentially contained in the works discussing the implications of averaging for asymptotic dynamics, see for example \cite{Hale,SV,GH}.

\subsection{Estimate for the variational problem}

In the context of the assumptions of
Lemma~\ref{lem:prob1-estm-norm} we will be interested  in an explicit
estimate for $\frac{\partial \varphi}{\partial x}(t_0,t,x) -
\frac{\partial\varphi_a}{\partial x}(t_0,t,x)$. We will use these bounds in Section~\ref{sec:srzed-eq}.

\begin{lemma}
\label{lem:prob1-estm-var}
The same assumptions and notation as in Lemma~\ref{lem:prob1-estm-norm}. Additionally we assume that for all $k \in J_v$ holds
\begin{eqnarray*}
  \sup_{s \in [t_0,t_0+t], z \in W} \left\|\frac{\partial^2 F }{\partial
  z^2}(s,z)\right\|&=& C(D_{zz}F) < \infty, \\
   \sup_{z,z_1 \in W, s \in [t_0,t_0+h]} \|D_z F(s,z) \frac{v_k}{\partial z}(s,z_1) \| &=& C\left(D_z F \frac{\partial v_k}{\partial
   z}\right) < \infty\\
   \sup_{z \in W, s \in [t_0,t_0+h]} \|D_z F(s,z)  \| &=& C(D_z F) < \infty\\
    \sup_{s \in [t_0,t_0+t], z \in W}  \left\|\frac{\partial^2 v_k}{\partial t \partial z}(s,z)\right\|
    &=&
   C\left(\frac{\partial^2 v_k}{\partial t \partial z} \right)< \infty \\
    \sup_{s \in [t_0,t_0+t], z \in W}  \left\|\frac{\partial^2 v_k}{\partial z^2}(s,z)\right\|
    &=&C\left(\frac{\partial^2 v_k}{ \partial z^2}  \right) < \infty.
\end{eqnarray*}

Let
\begin{eqnarray*}
  l_2&=& l + \sum_{k \in J_v} C(g_k) C\left( \frac{\partial v_k}{\partial z}
   \right) \\
   \omega &=&\inf_{k \in J_v} |\omega_k|,
\end{eqnarray*}
and let $\tilde{b}(t)=\sum_{k \in J_v} b_k(t)$ be as in the assertion of Lemma~\ref{lem:prob1-estm-norm}.

Then for $t \in [0,h]$ holds
\begin{equation}
   \left\| \frac{\partial \varphi}{\partial x} (t_0,t,x) - \frac{\partial \varphi_a}{\partial x} (t_0,t,x)   \right\|
     \leq \frac{B(t)}{\omega}
\end{equation}
where $B:[0,h] \to \mathbb{R}$ is continuous  function depending on $C(\dots)$ and constants $l$, $h$ given by
\begin{eqnarray*}
  B (t) \leq C(D_{zz}F) \tilde{b}(t) \left(\frac{e^{l_2 t} - e^{lt}}{l_2 - l} \right) +
    \sum_{k \in J_v} B_k(t),
\end{eqnarray*}
where
\begin{eqnarray*}
  B_k(t) =  C(G_k) \left( C\left(\frac{\partial v_k}{\partial z}\right)(e^{lt} + e^{l_2 t}) \right)  + \\
   C(G_k) \frac{e^{l_2 t}-e^{lt}}{l_2 - l} \cdot \left(  C\left(D_z F \frac{\partial v_k}{\partial z} \right)
      +  C\left(\frac{\partial v_k}{\partial z}\right)\left(
   C(D_z F) + \sum_{k \in J_v} C(g_i)C\left(\frac{\partial v_i}{\partial z}\right)
   \right) + \right. \\
     \left.
   +  \left( C\left(\frac{\partial^2 v_k}{\partial t \partial z} \right) +
  C\left(\frac{\partial^2 v_k}{ \partial z^2} \right) \cdot \left( C(F) + \sum_{i \in J_v} C(g_i) C(v_i)\right)
  \right)
      \right)
\end{eqnarray*}
\end{lemma}
\noindent
\textbf{Proof:}
Let us denote by $V_a(t,z)=\frac{\partial \varphi_a}{\partial
x}(t_0,t-t_0,z)$ and $V(t,z)=\frac{\partial \varphi}{\partial
x}(t_0,t-t_0,z)$.  Often in the remainder of the proof, we write
just $V_a(t)$ and $V(t)$ for $V_a(t,x(t_0))$ and $V(t,x(t_0))$. Observe that in this notation holds
\begin{eqnarray*}
  V_a(t_0)=V(t_0)=\mbox{Id}.
\end{eqnarray*}

 From (\ref{eq:auto-prob1-loc}) and
(\ref{eq:prob1-loc}) it follows that
\begin{equation}
\frac{dV_a}{dt}(t)=\frac{\partial F}{\partial y}(t,y(t)) V_a(t)
\label{eq:eqvar-va}
\end{equation}
and
\begin{equation}
\frac{dV}{dt}(t)=\frac{\partial F}{\partial y}(t,x(t)) V(t) +
\left(\sum_{k \in J_v} g_k(\omega_k t) \frac{\partial v_k}{\partial x}(t,x(t))\right)
V(t). \label{eq:eqvar-v}
\end{equation}

We rewrite (\ref{eq:eqvar-v}) as follows
\begin{eqnarray*}
\frac{dV}{dt}(t)=\frac{\partial F}{\partial y}(t,y(t)) V(t) +
\left(\frac{\partial F}{\partial y}(t,x(t)) -  \frac{\partial
F}{\partial y}(t,y(t)) \right) V(t)  + \\ \left(\sum_{k \in J_v} g_k(\omega_k t)
\frac{\partial v_k}{\partial x}(t,x(t))\right) V(t).
\end{eqnarray*}

We obtain the following equation for $\Delta(t)=V(t)-V_a(t)$
\begin{eqnarray*}
  \frac{d \Delta}{dt}= \frac{\partial F}{\partial y}(t,y(t)) \cdot
  \Delta + \left(\frac{\partial F}{\partial y}(t,x(t)) -  \frac{\partial
F}{\partial y}(t,y(t)) \right) V(t)  + \\ \left(\sum_{k \in J_v} g_k(\omega_k t)
\frac{\partial v_k}{\partial x}(t,x(t))\right) V(t).
\end{eqnarray*}
Let $M(t,t_0)$ be the fundamental matrix for (\ref{eq:eqvar-va}).
Observe that in our notation we have
$M(t+t_0,t_0)=V_a(t_0+t)=V_a(t+t_0,x(t_0))$.

We have
\begin{eqnarray*}
  \Delta(t+t_0) = S(t) + \sum_{k \in J_v} I_k(t),
\end{eqnarray*}
where
\begin{eqnarray*}
S(t)= \int_0^t M(t_0+t,t_0+s)  \left(\frac{\partial F}{\partial
y}(t_0+s,x(t_0+s)) -  \frac{\partial F}{\partial
y}(t_0+s,y(t_0+s)) \right) V(t_0+s) ds \\
I_k(t) = \int_0^t  g_k(\omega_k (t_0+s)) M(t_0+t,t_0+s) \frac{\partial
v_k}{\partial x}(t_0+s,x(t_0+s)) V(t_0+s) ds
\end{eqnarray*}

To provide bounds for $S(t)$ and $I_k(t)$ we need  estimates for
$V(t_0+s)$ and $M(t_0+t,t_0+s)$ for $s \in [0,t]$.
This can be obtained  using the logarithmic norm for
(\ref{eq:prob1-loc})  and (\ref{eq:auto-prob1-loc}), respectively (compare the proof of Lemma~\ref{lem:prob1-estm-var}).  We have
\begin{eqnarray*}
  \|V(t_0+s)\| \leq \exp(l_2 s)
\end{eqnarray*}
where
\begin{eqnarray*}
  \sup_{s \in [t_0,t_0+t], t_1 \in \mathbb{R}, x \in W}
  \left(\mu(D_xf(s,x)) + \sum_{k \in J_v} |g_k(t_1)| \left\| \frac{\partial v_k}{\partial x}(s,x) \right\|
  \right)\\
     \leq l + \sum_{k \in J_v} C(g_k) C\left( \frac{\partial v_k}{\partial z}
   \right)=:l_2
\end{eqnarray*}
In the above estimate obtained from (\ref{eq:prob1-loc}) we applied the logarithmic norm to $F(t,x)$ and the standard norm to the integral. This can be justified as follows.

We have $\mu(A+B)=2\mu\left((A+B)/2\right) \leq \mu(A) + \mu(B)$, the last inequality follows from the convexity of the logarithmic norm.

Since $\mu(B) \leq \|B\|$, so we obtain $\mu(A+B) \leq \mu(A) + \|B\|$.

The bound for $M(t_0+t,t_0+s)$ is given by
\begin{equation*}
\|M(t_0+t,t_0+s)\| \leq e^{l(t-s)}, \quad s \in [0,t].
\end{equation*}

We are now ready to estimate $S(t)$.

From Lemma~\ref{lem:prob1-estm-norm} it follows that for $\tilde{b}(t)=\max_{s \in [0,t]}\sum_{k \in J_v} b_k(s)$ holds
\begin{eqnarray*}
 \|S(t)\| \leq \int_0^t e^{l(t-s)}  C(D_{zz}F) \|x(t_0+s) - y(t_0+s)\| e^{l_2 s}
 ds \leq \\
 \frac{e^{lt} C(D_{zz}F) \tilde{b}(t)}{\omega} \int_0^t e^{(l_2 - l)s}  ds
 =
   \frac{C(D_{zz}F) \tilde{b}(t)}{\omega} \left(\frac{e^{l_2 t} - e^{lt}}{l_2 - l} \right)
\end{eqnarray*}

To estimate $I_k(t)$ we will integrate  by parts.

We have
\begin{eqnarray*}
 I_k = \left. \left(\frac{G_k(\omega_k(t_0+s))}{\omega_k}  M(t_0+t,t_0+s)  \frac{\partial
v_k}{\partial x}(t_0+s,x(t_0+s)) V(t_0+s)\right)
\right|_{s=0}^{s=t} + \\
- \frac{1}{\omega_k} \int_0^t G_k(\omega_k (t_0+s)) \frac{\partial
}{\partial s} \left( M(t_0+t,t_0+s) \frac{\partial v_k}{\partial
x}(t_0+s,x(t_0+s)) V(t_0+s) \right) ds
\end{eqnarray*}
To estimate  the second term we use the following identities (we use (\ref{eq:Mfund-4}))
\begin{eqnarray*}
  \frac{\partial
}{\partial s} \left( M(t_0+t,t_0+s) \frac{\partial v_k}{\partial
x}(t_0+s,x(t_0+s)) V(t_0+s) \right) = \\
-  M(t_0+t,t_0+s) \frac{\partial F}{\partial y}(t_0+s,y(t_0+s))
\frac{\partial v_k}{\partial
x}(t_0+s,x(t_0+s)) V(t_0+s) + \\
 M(t_0+t,t_0+s) \left(\frac{\partial^2 v_k}{\partial t\partial
x}(t_0+s,x(t_0+s)) +  \frac{\partial^2 v_k}{\partial
x^2}(t_0+s,x(t_0+s)) x'(t_0+s) \right) V(t_0+s)  + \\
  M(t_0+t,t_0+s) \frac{\partial v_k}{\partial x}(t_0+s,x(t_0+s))
V'(t_0+s).
\end{eqnarray*}
Therefore we obtain for $s \in [0,t]$
\begin{eqnarray*}
\left\| \frac{\partial }{\partial s} \left( M(t_0+t,t_0+s)
\frac{\partial v_k}{\partial x}(t_0+s,x(t_0+s)) V(t_0+s) \right)
\right\| \leq \\
 e^{l(t-s)} C\left(D_z F \frac{\partial v_k}{\partial z} \right) e^{l_2 s} + \\
 e^{l(t-s)}
  \left( C\left(\frac{\partial^2 v_k}{\partial t \partial z} \right) +
  C\left(\frac{\partial^2 v_k}{ \partial z^2} \right) \cdot \left( C(F) + \sum_{i \in J_v} C(g_i) C(v_i)\right)
  \right)
   e^{l_2 s} + \\
   e^{l(t-s)} C\left(\frac{\partial v_k}{\partial z}\right)\left(
   C(D_z F) + \sum_{i \in J_v} C(g_i)C\left(\frac{\partial v_i}{\partial z}\right)
   \right) e^{l_2 s}
\end{eqnarray*}
We obtain the following estimate for $I_k(t)$.
\begin{eqnarray*}
 |\omega_k| \cdot  \|I_k(t)\| \leq   C(G_k) \left( C\left(\frac{\partial v_k}{\partial z}\right)\left(e^{lt} + e^{l_2 t}\right) \right)  + \\
   C(G_k) \frac{e^{l_2 t}-e^{lt}}{l_2 - l} \cdot \left(  C\left(D_z F \frac{\partial v_k}{\partial z} \right)
      +  C\left(\frac{\partial v_k}{\partial z}\right)\left(
   C(D_z F) + \sum_{i \in J_v} C(g_i)C\left(\frac{\partial v_i}{\partial z}\right)
   \right) \right. +\\
     \left.
    + \left( C\left(\frac{\partial^2 v_k}{\partial t \partial z} \right) +
  C\left(\frac{\partial^2 v_k}{ \partial z^2} \right) \cdot \left( C(F) + \sum_{i \in J_v} C(g_i) C(v_i)\right)
  \right)
      \right)
\end{eqnarray*}
This finishes the proof. \qed

%In the above lemma we need an upper bound  on $\|D_z F(t,z)\|$,
%which in the context of the Galerkin projections  dissipative PDEs  of increasing dimension is unbounded. This appears because we need
%an estimate for $V'$. We believe it is possible to avoid this by
%a more careful analysis.

%\input hyperbolic.tex

%\input globattracting.tex

\section{Planar non-autonomous  equation with very strong expansion}
\label{sec:srzed-eq}

Let us write a real two-dimensional version of (\ref{eq:srzednicki})
\begin{eqnarray*}
  \frac{dx}{dt}&=& x(1+\cos (\kappa t) |z|^2) + \sin(\kappa t)|z|^2 y, \\
  \frac{dy}{dt}&=& -y(1+\cos (\kappa t) |z|^2) +  \sin(\kappa t)|z|^2 x.
\end{eqnarray*}
In this section we will use quite often the following notations $c=\cos(\kappa t)$, $s=\sin(\kappa t)$.
The norm on $\mathbb{C}=\mathbb{R}^2$ is the euclidian norm and it is used to define the operator norms for linear and bilinear map arising in the analysis of (\ref{eq:srzednicki}).

\subsection{All solutions with large initial data go  to infinity in finite time }

Consider first differential inequality
\begin{equation}
  \frac{dx}{dt} > b x^3.   \label{eq:1D-x3}
\end{equation}
\begin{lemma}
\label{lem:1D-explosion}
  Assume $b >0$ and $x:[t_0,t_{max})\to \mathbb{R}$ is $C^1$ function, which satisfies (\ref{eq:1D-x3}) and $x(t_0)=x_0 >0$. Then
  \begin{equation*}
    x(t) > \frac{x_0}{\sqrt{1 - 2bx^2_0(t-t_0)}}.
  \end{equation*}
  In particular, $t_{max}-t_0 \leq \frac{1}{2bx_0}$.
\end{lemma}
\noindent
\textbf{Proof:}
Let $y(t)$ be a solution of the following Cauchy problem
\begin{equation}
  \frac{dy}{dt} = b y^3, \quad y(t_0)=y_0.
\end{equation}
 It is easy to see that
 \begin{equation}
   y(t)=\frac{y_0}{\sqrt{1 - 2by_0^2 (t-t_0)}}
 \end{equation}

 From the standard differential inequalities it follows that if $x_0 \geq y_0$, then $x(t) > y(t)$ for $t > t_0$.
\qed

Now we will show the existence of the forward invariant cones for short times.  It will turn out that in such cone a solution might explode to infinity
in time shorter than the established time for the forward invariance.

We define a quadratic form $Q(x,y)=x^2 - y^2$. Let us define  cones $Q^+=\{z \ | \ Q(z) \geq  0\}$ and $Q^-=\{z \ | \ Q(z) \leq  0\}$.

\begin{lemma}
 \label{lem:srzed-q+finv}
  Cone $Q^+$ is forward invariant as long as
  \begin{eqnarray*}
       \cos(\kappa t)&>&0.
  \end{eqnarray*}
\end{lemma}
\noindent
\textbf{Proof:}
Let $z(t)=x(t)+iy(t)$ be a nonzero solution of (\ref{eq:srzednicki}).  We have
\begin{eqnarray*}
  \frac{1}{2}\frac{d}{dt} Q(x(t),y(t))=  xx' - yy' = \\
  x^2 (1+c|z|^2) + |z|^2 s xy + y^2 (1+c|z|^2) -  |z|^2 s xy = |z|^2 (1+c|z|^2) \geq c |z|^4
\end{eqnarray*}
\qed

\begin{lemma}
\label{lem:x'-in-pos-cone}
Assume that $z=(x+iy) \in Q^+$, $x>0$ and $\cos(\kappa t) - |\sin(\kappa t)|>0$, then
\begin{equation*}
  \frac{dx}{dt} >  x^3 (\cos(\kappa t) -  |\sin(\kappa t)|).
\end{equation*}
\end{lemma}
\noindent
\textbf{Proof:}
Since $|z|^2 \geq x^2$ and $|y| \leq  x$  we have
\begin{eqnarray*}
  \frac{dx}{dt} > c |z|^2 x - |s| |z|^2 |y| \geq  |z|^2 (cx - |s|x)= x |z|^2 (c - |s|) \geq x^3 (c-|s|)
\end{eqnarray*}
\qed

\begin{lemma}
\label{lem:toInfinityInFiniteTime}
Let $t_1 >0$ and $\epsilon >0$ be such that the following inequality is satisfied for $t \in [0,t_1]$
\begin{eqnarray*}
\cos(\kappa t) - |\sin(\kappa t)|> \epsilon > 0.
\end{eqnarray*}

Assume that $z_0=(x_0,y_0) \in Q^+$ and $x_0^2 \geq \frac{1}{2 \epsilon t_1}$. Then, $z(t)$, the solution of (\ref{eq:srzednicki}) with the initial condition $z(0) = z_0$
goes to infinity in finite time.
\end{lemma}
\noindent
\textbf{Proof:}
Without any loss of the generality we can assume that $x_0>0$.

From Lemmas~\ref{lem:srzed-q+finv} and \ref{lem:x'-in-pos-cone} it follows that for $t \in [0,t_1]$ holds
\begin{equation}
  \frac{dx}{dt} > \epsilon x^3.
\end{equation}
From Lemma~\ref{lem:1D-explosion} if follows that for $t \in [0,t_1]$ holds
\begin{equation}
 x(t) > \frac{x_0}{\sqrt{1 - 2\epsilon x^2_0 t}}.
\end{equation}
Observe that we have the solution blow-up before time $t_1$ if
\begin{equation}
   x_0^2 \geq \frac{1}{2 \epsilon t_1}.
\end{equation}

\qed

\begin{theorem}
\label{thm:toInfinityInFiniteTime}
Let $\delta \in (0,\pi/4)$. If $|z_0|^2 \geq \frac{\kappa}{ \sqrt{2} (\pi/4 - \delta) \sin \delta}$, then for any $t_0$ the solution of (\ref{eq:srzednicki})
with initial condition $z(t_0)=z_0=(x_0,y_0)$ goes to infinity in finite time  backwards or forwards in time.
\end{theorem}
%\textbf{Zobaczyc, gdzie maksymalizuje sie mianownik w oszacowaniu na $|z_0|^2$}
\noindent
\textbf{Proof:}
Let us fix  $t_0=0$. Let $t_1=\frac{\pi/4 - \delta}{\kappa}$
Then for $t \in [0,t_1]$ holds
\begin{eqnarray*}
  \cos (\kappa t) -  |\sin(\kappa t)|  =  \cos (\kappa t_1) - \sin(\kappa t_1) \geq \\
  \cos (\pi/4 - \delta) - \sin(\pi/4 - \delta)= \sqrt{2} \sin \delta.
\end{eqnarray*}
From Lemma~\ref{lem:toInfinityInFiniteTime} with $\epsilon=  \sqrt{2} \sin \delta $ we obtain that if the initial condition $z_0=(x_0,y_0)$ satisfies
\begin{equation}
  x_0^2 \geq \frac{1}{2 \epsilon t_1}= \frac{\kappa}{2 \sqrt{2} (\pi/4 - \delta) \sin \delta}
\end{equation}
then $|x(t)| \to \infty$ in finite time.

Observe that by inverting time we obtain that for initial conditions in $Q_1^- $  with $y_0^2 \geq  \frac{\kappa}{2 \sqrt{2} (\pi/4 - \delta) \sin \delta}$  the solutions will explode backward in time. Therefore we see that if $x_0^2 + y_0^2 \geq \frac{\kappa}{ \sqrt{2} (\pi/4 - \delta) \sin \delta}$, then
we have  $|x(t)| \to \infty$ forward in finite time or $|y(t)| \to \infty$ backward in finite time.

In the above reasoning $t_0=0$ was singled out. But in fact there is nothing special about this initial time. The transformation
\begin{equation}
  z_1=\exp(-i\kappa t_0/2) z
\end{equation}
transforms (\ref{eq:srzednicki}) into
\begin{equation}
  z_1'=\overline{z}_1 (1 + |z_1|^2 \exp(i \kappa (t-t_0))).
\end{equation}

\qed

\subsection{A priori estimates for small time step}

\begin{lemma}
\label{lem:srze-rough-encl}
  Let $\varphi$ be the local process induced by (\ref{eq:srzednicki}) and $R_0 \geq 1$. For any $t_0 \in \mathbb{R}$ if $|z_0| \leq R_0$ and  $h=\frac{1}{8R_0^2}$, then
  $\varphi(t_0,t,z_0)$ is defined for $t \in [0,h]$ and
  \begin{equation}
    |\varphi(t_0,t,z_0)| \leq \sqrt{2} R_0.
  \end{equation}
\end{lemma}
\noindent \textbf{Proof:}
It is easy to see that for $|z| \geq 1$ holds
\begin{equation}
   \frac{d |z|}{dt} \leq |z| (1 + |z|^2 ) \leq 2 |z|^3.
\end{equation}
From Lemma~\ref{lem:1D-explosion} with $b=2$ we obtain the following estimate for $|z_0| \leq R_0 $ and $t \in [0,h]$
\begin{equation}
  |\varphi(t_0,t,z_0)| \leq \frac{R_0}{\sqrt{1 - 2 b R_0^2 t}}  \leq \sqrt{2} R_0
\end{equation}
\qed

\comment{
\textbf{Nie wiem czy to zatrzymac:}
We could do a bit better estimate in the above lemma as follows:
\begin{eqnarray*}
  \frac{d|z|}{dt} \leq |z| + |z|^3.
\end{eqnarray*}
Since  solutions of equation
\begin{equation}
  \frac{dx}{dt}=x+x^3,
\end{equation}
are given by
\begin{equation}
  x(t)=\frac{k\exp(t-t_0)}{\sqrt{1 - k^2 \exp(2(t-t_0))}}, \quad k=\frac{x(t_0)}{\sqrt{1+x(t_0)^2}}
\end{equation}
\textbf{Troche to ciezko ladnie szacowac, $k \to 1$ odpowiada za wzrost $x_0$ }
}

\subsection{$C^0$ estimates for the averaging}
We would like to write down explicitly the estimates from Lemma~\ref{lem:prob1-estm-norm} for (\ref{eq:srzednicki}).

Let us first rewrite equation (\ref{eq:srzednicki}) in the form (\ref{eq:prob1-loc}).  We have
\begin{equation}
  z'=f(z) + g_1(\omega_1 t) v_1(z) + g_2(\omega_2 t) v_2(z),
\end{equation}
where
\begin{eqnarray*}
  f(z)&=&\overline{z}, \\
  \omega_1&=&\omega_2=\kappa, \\
  g_1(t)&=&\cos t, \quad G_1(t) =  \sin t \\
  g_2(t)&=&\sin t, \quad G_2(t) = - \cos t \\
  v_1(z)&=&|z|^2 (x,-y)^t, \\
  v_2(z)&=&|z|^2 (y,x)^t.
\end{eqnarray*}
Therefore we see that $J_v=\{1,2\}$ and on $J_v$ the counting measure is used.

We will also need the estimates for the partial derivatives of $v_k$.
\begin{equation*}
  \frac{\partial v_1}{\partial z}=\left[\begin{array}{cc}
                                    3x^2 + y^2, & 2xy \\
                                    -2xy, & -x^2 - 3y^2
                                  \end{array}\right],
   \qquad
        \frac{\partial v_2}{\partial z}=\left[\begin{array}{cc}
                                    2xy, & x^2 + 3y^2 \\
                                    3x^2 + y^2, & 2xy
                                  \end{array}\right].
\end{equation*}

\begin{lemma}
\label{lem:estmDvi}
For $z \in \mathbb{R}^2$ holds
\begin{equation}
  \|D_z v_1(z)\|=\|D_z v_2 (z)\|=3r^2
\end{equation}
\end{lemma}
\noindent
\textbf{Proof:}
Let $z=(x,y)=r(\cos \alpha, \sin\alpha)$. It is easy to see that
\begin{eqnarray*}
  D v_1(z) &=& \left[\begin{array}{cc}
                                    2r^2 + (x^2 -y^2), & 2xy \\
                                    -2xy, & -2r^2 +  (x^2 - y^2)
                                  \end{array}\right]= \\
                                  &=& r^2 \left(
                                   \left[\begin{array}{cc}
                                    2 , & 0 \\
                                    0, & -2
                                  \end{array}\right] +
                                   \left[\begin{array}{cc}
                                    \cos (2\alpha), & \sin(2\alpha) \\
                                    -\sin(2\alpha), & \cos(2\alpha)
                                  \end{array}\right] \right),
\end{eqnarray*}
and
\begin{eqnarray*}
  D v_2(z) &=& \left[\begin{array}{cc}
                                    2xy, & 2r^2 + (y^2 - x^2) \\
                                    2r^2 + x^2 - y^2, & 2xy
                                  \end{array}\right]= \\
                                  &=& r^2 \left(
                                   \left[\begin{array}{cc}
                                    0 , & 2 \\
                                    2, & 0
                                  \end{array}\right] +
                                   \left[\begin{array}{cc}
                                    \sin (2\alpha), & -\cos(2\alpha) \\
                                    \cos(2\alpha), & \sin(2\alpha)
                                  \end{array}\right] \right),
\end{eqnarray*}
Therefore we have
\begin{equation*}
  Dv_i(z) = r^2 (A_i + Q_i), \quad i=1,2
\end{equation*}
where $A_i$ is a symmetric matrix with eigenvalues $\pm 2$ and $Q_i$ is a rotation matrix.  Therefore we have
\begin{equation}
  \|Dv_i(z)\|  \leq r^2 (\|A_i\| + \|Q_i\| )=3r^2.
\end{equation}
The proof of equality is left to the reader.
\qed

Therefore we have the following bounds on $W=\overline{B}(0,R)$ and $t \in [t_0,t_0+h]$, where $h=\frac{1}{4R^2}$ (compare with Lemma~\ref{lem:srze-rough-encl},
there $R_0$ was the size of the initial condition and here $R=\sqrt{2} R_0$ is the size of the enclosure)
\begin{eqnarray*}
  C(f)&=&R \\
  C(G_k)&=&C(g_k)=1, \quad k=1,2 \\
  l&=&1, \\
  C(v_k)&=& C(D_z f v_k)= R^3,  \quad k=1,2\\
  C\left( \frac{\partial v_k}{\partial t}\right) &=&0, \quad k=1,2 \\
  C\left(\frac{\partial v_k}{\partial z} \right) &\leq& M R^2,  \quad k=1,2, \\
  C\left(D_z v_k \tilde{F} \right)&=&  C\left(\frac{\partial v_k}{\partial z} \right) \left( C(f)+\sum_{i \in J_v} C(v_i)C(g_i)\right)
\end{eqnarray*}
From Lemma~\ref{lem:estmDvi} it follows that $M=3$.

Hence for $k=1,2$ from Lemma~\ref{lem:prob1-estm-norm} we obtain
\begin{eqnarray*}
   b_k(t) \leq  C(G_k) \left( C(v_k) (1+e^{lt})   +  \right.\\
  \left.  \left( C\left(D_z fv_k\right) + C\left(\frac{\partial v_k}{\partial t}\right) + C\left(\frac{\partial v_k}{\partial z}\right)\left(C(f)+\sum_{i \in J_v} C(v_i)C(g_i)\right) \right)
\frac{e^{lt}-1}{l} \right) = \\
  R^3 (1+e^{t})   +  ( R^3  + M R^2 (R+ 2R^3))(e^{t}-1) = R^3((2+M)e^t - M) + 2MR^5(e^t - 1)
\end{eqnarray*}
Summarizing, we have proved
\begin{lemma}
\label{lem:eq-srzed-c0-estm}
For $t \in \left[0,\frac{1}{4R^2} \right]$ and $|z_0| \leq \frac{R}{\sqrt{2}}$ holds
\begin{eqnarray*}
  |\varphi(t_0,t,z_0) - (e^t x_0,e^{-t}y_0)| \leq \frac{\tilde{b}(t)}{|\kappa|}:=\frac{2}{|\kappa|} \left(  R^3((2+M)e^t - M) + 2MR^5(e^t - 1) \right).
\end{eqnarray*}
\end{lemma}

\subsection{$C^1$ estimates}

We want to use Lemma~\ref{lem:prob1-estm-var}. For this we need to estimate also the second derivatives of $v_k$.

We have
\begin{eqnarray*}
   \frac{\partial^2 v_{1,x}}{\partial z^2}&=& \left[\begin{array}{cc}
                                    6x, & 2y \\
                                    2y, & 2x
                                  \end{array}\right],
          \qquad   \frac{\partial^2 v_{1,y}}{\partial z^2}= \left[\begin{array}{cc}
                                    -2y, & -2x \\
                                    -2x, & -6y
                                  \end{array}\right],         \\
       \frac{\partial^2 v_{2,x}}{\partial z^2}&=& \left[\begin{array}{cc}
                                    2y, & 2x \\
                                    2x, & 6y
                                  \end{array}\right],
          \qquad   \frac{\partial^2 v_{2,y}}{\partial z^2}= \left[\begin{array}{cc}
                                    6x, & 2y \\
                                    2y, & 2x
                                  \end{array}\right].
\end{eqnarray*}
To estimate $ C\left(\frac{\partial^2 v_k}{ \partial z^2}  \right)$ we will use the following lemma.

\begin{lemma}
\label{lem:eq-srz-d2v-estm}
  For $z \in \mathbb{R}^2$ and $k=1,2$ it  holds
  \begin{equation}
    \left\|  \frac{\partial^2 v_{k}}{\partial z^2} \right\| \leq r \sqrt{24 + 16\sqrt{2}}= \left(4 + 2 \sqrt{2}\right)r,
  \end{equation}
  where $r=|z|$,
\end{lemma}
\textbf{Proof:}
Fist we will  estimate
\begin{equation}
 C_{1,x}=\left\|  \frac{\partial^2 v_{1,x}}{\partial z^2} \right\|
 \end{equation}
For this  it is enough to estimate from above the spectrum of the matrix $\left[\begin{array}{cc}
                                    6x, & 2y \\
                                    2y, & 2x
                                  \end{array}\right]$.
 An easy computation yields the following formula for the eigenvalues
 \begin{equation*}
   \lambda= 4x \pm  2r.
 \end{equation*}
 Therefore we obtain
 \begin{equation*}
   C_{1,x} = 4 |x| + 2r \leq 6r.
 \end{equation*}
 In order to estimate $D^2 v_1$ observe first that from symmetry arguments it follows immediately that
 \begin{equation*}
   C_{1,y}=\left\|  \frac{\partial^2 v_{1,y}}{\partial z^2} \right\|= 4|y| +  2r.
 \end{equation*}
 Therefore we have for any vectors $a,b \in \mathbb{R}^2$
 \begin{eqnarray*}
   \|D^2 v_1 (a,b)\|^2 \leq (C_{1,x} \|a\| \cdot  \|b\|)^2 +   (C_{1,y} \|a\| \cdot  \|b\|)^2 = \left( C_{1,x}^2 + C_{1,y}^2 \right) (\|a\| \cdot  \|b\|)^2 \\
   \|D^2 v_1\| \leq \left( C_{1,x}^2 + C_{1,y}^2 \right)^{1/2}.
 \end{eqnarray*}
 We have
 \begin{eqnarray*}
    C_{1,x}^2 + C_{1,y}^2 = (16 x^2 + 16 |x| r  + 4 r^2 ) + (16 y^2 + 16 |y| r  + 4 r^2 )=\\
    24r^2 + 16(|x|+|y|)r \leq 24r^2 + 16 \sqrt{2} r^2=
        (24 + 16\sqrt{2})r^2.
 \end{eqnarray*}

\qed

Let $h$, $W$ be as in the previous subsection. We have
\begin{eqnarray*}
   C(D_{zz}f) &=& 0, \\
   C\left(D_z f \frac{\partial v_k}{\partial z}\right) &=& C\left(\frac{\partial v_k}{\partial z}\right) = MR^2,\\
    C(D_z f) &=& 1, \\
   C\left(\frac{\partial^2 v_k}{\partial t \partial z} \right) &=&  0, \\
  C\left(\frac{\partial^2 v_k}{ \partial z^2}  \right) &=& NR.
\end{eqnarray*}
From Lemma~\ref{lem:eq-srz-d2v-estm} it follows that $N=4+2 \sqrt{2}$.

We have
\begin{eqnarray*}
  l_2&=& l + \sum_{k \in J_v} C(g_k) C\left( \frac{\partial v_k}{\partial z}
   \right) = 1+ 2 MR^2\\
   \omega &=&|\kappa|,
\end{eqnarray*}
and  $\tilde{b}(t)$  given in Lemma~\ref{lem:eq-srzed-c0-estm}.

From Lemma~\ref{lem:prob1-estm-var} we obtain (we set $k=1$ below)
\begin{eqnarray*}
  B(t)=2B_1(t)= 2 C(G_k)\left(C\left(\frac{\partial v_k}{\partial z}\right)(e^{lt}+e^{l_2t}) \right)+ \\
  2 C(G_k) \frac{e^{l_2 t}-e^{lt}}{l_2 - l} \cdot \left(  C\left(D_z f \frac{\partial v_k}{\partial z} \right)
      +  C\left(\frac{\partial v_k}{\partial z}\right)\left(
   C(D_z f) + \sum_{i \in J_v} C(g_i)C\left(\frac{\partial v_i}{\partial z}\right)
   \right) \right. + \\
    + \left.
      \left( C\left(\frac{\partial^2 v_k}{\partial t \partial z} \right) +
  C\left(\frac{\partial^2 v_k}{ \partial z^2} \right) \cdot \left( C(f) + \sum_{i \in J_v} C(g_i) C(v_i)\right)
  \right)
      \right) =
 2 MR^2 \left(e^t + e^{(1+2MR^2)t}\right)+  \\2 \frac{e^{(1+2MR^2) t}-e^{t}}{2MR^2} \cdot \left(  MR^2   + MR^2( 1 + 2 MR^2) +  NR  \cdot ( R  + 2 R^3) \right) = \\
 2 MR^2 e^t\left( 1 + e^{2MR^2 t}\right)+  e^t(e^{2M R^2 t }-1)    \left( 2 + \frac{N}{M} + 2 R^2 \left(M + \frac{N}{M}\right) \right)
\end{eqnarray*}
We proved the following
\begin{lemma}
\label{lem:eq-srzed-c1-estm}
For $t \in \left[0,\frac{1}{4R^2} \right]$ and $|z_0| \leq \frac{R}{\sqrt{2}}$ holds
\begin{eqnarray*}
  \left\|\frac{\partial \varphi}{\partial z_0}(t_0,t,z_0) -
  \left[\begin{array}{cc}
    e^t & 0 \\
    0 & e^{-t}
  \end{array}
  \right] \right\| \leq
    \frac{1}{|\kappa|}  2 MR^2 e^t\left( 1 + e^{2MR^2 t}\right) + \\
     \frac{e^t(e^{2MR^2 t }-1)}{|\kappa|}    \left( 2 + \frac{N}{M} + 2R^2 \left(M + \frac{N}{M}\right) \right)=: \frac{\tilde{B}(t,R)}{|\kappa|}
\end{eqnarray*}
\end{lemma}

\subsection{The estimates for the domain of hyperbolicity}

Let
\begin{equation*}
A= \left[\begin{array}{cc}
    1 & 0 \\
    0 & -1
  \end{array}
  \right].
\end{equation*}
For any $h >0$ map $e^{Ah}$ is hyperbolic and represents the linearization of (\ref{eq:srzednicki}) at $0$. Now we add to $e^{Ah}$ the influence of oscillating
terms from (\ref{eq:srzednicki}) and we would like to find possibly large $R$, such that in the ball of radius $R$  the shift by $h$ along trajectory of (\ref{eq:srzednicki}) is a hyperbolic map, in the sense discussed below.

\subsubsection{Cone conditions}

Let us define two cone fields, for $z\in \mathbb{R}^2$ we set $Q^+(z)=\{ z+(x,y) | |x| \geq |y|\}$, $Q^-(z)=\{ z+(x,y) | |x| \leq |y|\}$.

\begin{definition}
\label{def:hyper}
Let $\{F_j\}_{j \in \mathbb{Z}}$ be a family of maps,
$F_j:\mathbb{R}^2 \supset \mbox{dom} F_j \to \mathbb{R}^2$ be a $C^1$-map, such $F_j(0)=0$.  Let $N \subset \mbox{dom} F_j$ (for $j \in \mathbb{Z}$)  be a connected open set, such that $0 \in N$. We will say that family of maps $\{F_j\}$ is hyperbolic on $N$ iff the following conditions hold:
\begin{itemize}
\item the cone field $Q^+$ is forward invariant relatively to $N$, i.e.
  if $z_1,z_2 \in N$ are such that $z_2 \in Q^+(z_1)$, then for all $j \in \mathbb{Z}$ $F_j(z_2) \in Q^+(F_j(z_1))$
\item there exist constants $0 \leq \mu < 1 < \xi$, such that
\begin{itemize}
\item if $z_1, z_2 \in N$ and $z_2 \in Q^-(z_1)$, then for all $j \in \mathbb{Z}$
   \begin{equation}
     |\pi_y (F_j(z_1) - F_j(z_2))|  \leq  \mu |\pi_y (z_1 - z_2)|, \label{eq:contr-neg-cone}
   \end{equation}
\item    if $z_1, z_2 \in N$ and $z_2 \in Q^+(z_1)$, then  for all $j \in \mathbb{Z}$
   \begin{equation}
     |\pi_x (F_j(z_1) - F_(z_2))|  \geq  \xi |\pi_x (z_1 - z_2)|. \label{eq:exp-pos-cone}
   \end{equation}
\end{itemize}
\end{itemize}
\end{definition}

We have the following simple theorem, which is an adaptation of results from \cite{CaZ}, where normally hyperbolic invariant manifolds are discussed, (see also \cite{ZCC})
\begin{theorem}
\label{thm:hyper}
Assume that family of maps $\{F_j\}$  is a hyperbolic on a convex set $N$ in the sense of Definition~\ref{def:hyper} and assume that for all $j \in \mathbb{Z}$ $F_j^{-1}(0)=\{0\}$. Then
\begin{itemize}
\item if $z \in N \setminus \{0\}$, $l \in \mathbb{Z}$ and $z \in Q^+(0)$, then there exists $n \geq 1$, such that
\begin{eqnarray*}
   F_{k+l} \circ F_{k-1+l}\circ \cdots \circ F_l (z) &\in& Q^+(0) \cap N, \quad \mbox{for} \quad  k=1,\dots,n-1, \\
   F_{n+l} \circ F_{n-1+l}\circ \cdots \circ F_l (z) &\notin& N.
\end{eqnarray*}
\item if $z \in  N \setminus \{0\}$, $l \in \mathbb{Z}$, and $z \in Q^-(0)$, then one of the two following conditions is satisfied
   \begin{itemize}
     \item there exists $n \geq 1$   such that
       \begin{eqnarray*}
      F_{k+l} \circ F_{k-1+l}\circ \cdots  \circ F_l(z) &\in& N \cap Q^-(0), \mbox{ for $k=1,\dots,n-1$} \\
       F_{n+l} \circ F_{n-1+l}\circ \cdots \circ F_l (z) &\in& Q^+(0).
       \end{eqnarray*}
     \item
       \begin{eqnarray*}
       F_{k+l} \circ F_{k-1+l}\circ \cdots \circ F_l (z) &\in& N \cap Q^-(0), \quad \forall k \in \mathbb{N}\\
       \lim_{k \to \infty} F_{k+l} \circ F_{k-1+l}\circ \cdots \circ F_l(z)&=&0
       \end{eqnarray*}
   \end{itemize}
\end{itemize}

In particular, if $N$ is bounded, then for every point $z \neq 0$
 every maximal orbit of $\{F_j\}$ through $z$ leaves $N$ either forward or backward in time.
\end{theorem}

Below we give the criterion for the  hyperbolicity of family $\{F_j\}_{j \in \mathbb{Z}}$.
\begin{lemma}
\label{lem:cc-hyper}
Let for $j \in \mathbb{Z}$  $F_j:\mathbb{R}^2 \supset \mbox{dom} F_j \to \mathbb{R}^2$ be a $C^1$-map, such $F_j(0)=0$.  Let $\overline{N} \subset \mbox{dom} F_j$ for all $j \in \mathbb{Z}$ and $N$ be a convex bounded open set, such that $0 \in N$.

Let us set
\begin{eqnarray}
  \xi&=&\inf_{j \in \mathbb{Z}}\left( \inf_{z \in N} \left|\frac{\partial F_x}{\partial x}(z)\right| - \sup_{z \in N} \left|\frac{\partial F_x}{\partial y}\right|\right), \\
  \mu&=&\sup_{j \in \mathbb{Z}}\left(\sup_{z \in N} \left|\frac{\partial F_y}{\partial y}(z)\right| + \sup_{z \in N} \left|\frac{\partial F_y}{\partial x}\right|\right).
\end{eqnarray}
If $\xi > 1 > \mu$, then $\{F_j\}_{j \in \mathbb{Z}}$ is hyperbolic on $N$.
\end{lemma}
\textbf{Proof:}
We will present the proof for $F_j=F$. The generalization to a family of maps is trivial.

Let $z_1=(x_1,y_1)$ and $z_0=(x_0,y_0)$. Then we have
\begin{eqnarray}
  F_x(z_1) - F_x(z_0) &=& \int_0^1 \frac{\partial F_x}{\partial x}(z_0+t(z_1-z_0))dt \cdot (x_1-x_0) + \nonumber \\
   & & \int_0^1 \frac{\partial F_x}{\partial y}(z_0+t(z_1-z_0))dt \cdot (y_1-y_0), \label{eq:diffFx} \\
  F_y(z_1) - F_y(z_0) &=& \int_0^1 \frac{\partial F_y}{\partial x}(z_0+t(z_1-z_0))dt \cdot (x_1-x_0) + \nonumber \\
   & & \int_0^1 \frac{\partial F_y}{\partial y}(z_0+t(z_1-z_0))dt \cdot (y_1-y_0) \label{eq:diffFy}.
\end{eqnarray}

From conditions (\ref{eq:diffFx}), (\ref{eq:diffFy}) and the convexity of $N$ it follows that
\begin{eqnarray}
  |F_x(z_1) - F_x(z_0)| &\geq& \inf_{z \in N} \left|\frac{\partial F_x}{\partial x}(z)\right|  \cdot |x_1-x_0| - \sup_{z \in N} \left|\frac{\partial F_x}{\partial y}(z)\right| \cdot |y_1-y_0|  \\
  |F_y(z_1) - F_y(z_0)| &\leq& \sup_{z\in N} \left|\frac{\partial F_y}{\partial x}(z)\right| \cdot |x_1-x_0| + \sup_{z\in N} \left|\frac{\partial F_y}{\partial y}(z)\right| \cdot |y_1-y_0|. \label{eq:diff-Fy-1}
\end{eqnarray}
Therefore for $z_1 \in Q^+(z_0)$ (recall that $|x_1-x_0| \geq |y_1-y_0|$) we obtain
 \begin{eqnarray*}
  |F_x(z_1) - F_x(z_0)| &\geq& \xi |x_1 - x_0| \\
  |F_y(z_1) - F_y(z_0)| &\leq&  \mu |x_1-x_0|.
\end{eqnarray*}
We see that in this case
\begin{equation}
  |F_x(z_1) - F_x(z_0)| \geq \xi |x_1 - x_0| \geq \mu |x_1-x_0| \geq  |F_y(z_1) - F_y(z_0)|.
\end{equation}
Hence $F(z_1) \in Q^+(F(z_0))$. This establishes the forward cone invariance of the cone field $Q^+$ and  expansion condition (\ref{eq:exp-pos-cone}) in positive cone.

It remains to show  (\ref{eq:contr-neg-cone}) for $z_1 \in Q^-(z_0)$. From (\ref{eq:diff-Fy-1}) it follows that
\begin{equation*}
   |F_y(z_1) - F_y(z_0)| \leq  \mu |y_1-y_0|.
\end{equation*}

\qed

\subsubsection{Estimation of cone conditions for (\ref{eq:srzednicki})}

Let us fix any $t_0 \in \mathbb{R}$.
We will show that the family of maps $\{z \mapsto \varphi(t_0+jh,h,z)\}_{j \in \mathbb{Z}}$ is hyperbolic on some
ball $B(R)$.

Let us set
\begin{equation}
  \Delta(t_0,h,z)=\frac{\partial \varphi}{\partial z}(t_0,h,z) - e^{Ah}
\end{equation}

Observe that from Lemma~\ref{lem:eq-srzed-c1-estm} we have the following bounds, which do not depend on $t_0$, on $\|\Delta(t_0,h,z)\|$ for $|z| \leq \frac{R}{\sqrt{2}} $ and $h \in \left(0,\frac{1}{4R^2}\right]$ 
\begin{equation}
  \|\Delta(t_0,h,z)\| \leq \frac{\tilde{B}(h,R)}{|\kappa|}.
\end{equation}

Let us now compute $\xi$ and $\mu$ from Lemma~\ref{lem:cc-hyper}.  We have
\begin{eqnarray*}
  \xi \geq e^h - \frac{2\tilde{B}(h,R)}{|\kappa|}, \\
  \mu \leq e^{-h} + \frac{2\tilde{B}(h,R)}{|\kappa|}.
\end{eqnarray*}
Therefore for the hyperbolicity we need to satisfy the following conditions
\begin{eqnarray}
 \frac{2\tilde{B}(h,R)}{|\kappa|} < e^h -1 \\
 \frac{2\tilde{B}(h,R)}{|\kappa|} < 1 -e^{-h}.
\end{eqnarray}
Since for $h>0$ holds
\begin{equation*}
  e^h - 1 > 1 -e^{-h},
\end{equation*}
hence we are left with the following condition
\begin{equation}
  \frac{2\tilde{B}(h,R)}{|\kappa|} < 1 -e^{-h},
\end{equation}
hence
\begin{equation}
  |\kappa| > \frac{2 \tilde{B}(h,R)}{1-e^{-h}}. \label{eq:kappa-norm}
\end{equation}

In Table~\ref{tab:kappa-od-r} we list several values of $\kappa$ depending on $R_0=\frac{R}{\sqrt{2}} $
\begin{table}[h!]
\begin{equation*}
        \begin{array}{|c|c|c|}
        \hline
         R_0  & \kappa\\
        \hline\hline
         1 & 3655 \\
         10 & 2.24 \cdot 10^7 \\
        100 & 2.23 \cdot 10^{11}\\
        \hline
        \end{array}
    \end{equation*}
    \caption{The values of $\kappa$, such that in $B(0,R_0)$ the behavior of (\ref{eq:srzednicki}) is hyperbolic. $\kappa$ is computed from (\ref{eq:kappa-norm}) with $R=\sqrt{2}R_0$ from Lemma~\ref{lem:eq-srzed-c1-estm}.  We used $M=3$, $N=4+2\sqrt{2}$, $h=\frac{1}{4R^2}$.}
    \label{tab:kappa-od-r}
\end{table}

Therefore we have proved the following result.
\begin{theorem}
\label{thm:eq-srzed-hyperbolic}
For any $R>0$ let $R_0=\frac{R}{\sqrt{2}}$, $h=\frac{1}{4R^2_0}$  and assume that
\begin{equation}
  |\kappa| >  \frac{2 \tilde{B}(h,R)}{1-e^{-h}}.
\end{equation}
Then  for any $t_0$ the family of maps $\{z \mapsto \varphi(t_0+jh,h,z)\}_{j \in \mathbb{Z}}$, is hyperbolic on $B(0,R_0)$.
\end{theorem}

Let us estimate now the growth of $\kappa$, such that we have the hyperbolic behavior on $B(0,R_0)$.  It is easy to see that for $h=\frac{1}{8R_0^2}$
\begin{eqnarray*}
  \tilde{B}(h(R_0),R_0\sqrt{2}) &=& O(R_0^2) \\
  |\kappa| &\geq& O(R_0^2)/(1-e^{-h}) = \frac{O(R_0^2)}{\frac{1}{8R_0^2}}= O(R_0^4)
\end{eqnarray*}

\begin{corr}
\label{corr:lowerbnd}
There exists $K >0$, such that
if $z(t)$ is bounded orbit for (\ref{eq:srzednicki}),  then either $z(t)\equiv 0$ or there exists $t_1$, such that $|z(t_1)| \geq K \kappa^{1/4}$.
\end{corr}

\subsection{Bounds on the nontrivial bounded solutions }
\label{subsec:srzeeq-bnds}
In the investigations of (\ref{eq:srzednicki}) in \cite{SW,WZ,OW} the chaotic behavior was obtained by the constructions of suitably matched isolating segments,
called $(U,U^-)$ and $(W,W^-)$ in \cite{SW,WZ}, which satisfy $U \subset W$.
Segment $W$ exists for any $\kappa \neq 0$ and its existence implies the existence of nonzero periodic orbit in it.
The shape of $W$ in the extended phasespace is given by a square $[-R,R]^2$ rotating with angular velocity $\kappa/2$ (see Section 5 in \cite{WZ}),
where $R$ satisfies the following inequality \cite[Lem. 13]{WZ}
\begin{equation*}
  R^2(R^2-1) > (1/2 + \kappa/4)^2.
\end{equation*}
Therefore $2 \pi/\kappa$ periodic orbit related to this isolating segment exists for $R = O(\kappa^{1/2})$. This is in fact the bound  for any bounded orbit in this segment.

Combining this with Cor.~\ref{corr:lowerbnd} we obtain these periodic have to satisfy
\begin{eqnarray*}
  \forall t \ |z(t)| \leq O(\kappa^{1/2}) \\
  \exists t_1 \ |z(t_1)| \geq O(\kappa^{1/4}).
\end{eqnarray*}

Observe that from Theorem~\ref{thm:toInfinityInFiniteTime} we obtain also that all orbits such that $|z(t_0)| > O(\kappa^{1/2})$ escape to infinity forward
or backward in time.  Therefore there is not point in taking substantially larger $R$ for the isolating block.

From the above observations it follows that any bounded nonzero orbit of (\ref{eq:srzednicki}) has to satisfy the following conditions
\begin{eqnarray*}
  \forall t \ |z(t)| \leq O(\kappa^{1/2}) \\
  \exists \{t_n\}_{n \in \mathbb{Z}}, \lim_{n \to -\infty} t_n=-\infty, \lim_{n \to \infty} t_n=\infty \  \forall n \quad |z(t_n)| \geq O(\kappa^{1/4}).
\end{eqnarray*}

%\textbf{Pytanie: czy nie da sie rozszerzyc jakos znaczaco obszaru hyperbolicznosci poprzez optymalizacje wyboru $h(R)$. Mozliwe ze tak, ale raczej przez %rozwazanie przesuniecia o okres $\frac{2\pi}{\kappa}$. W tym celu trzeba przerobic lematy o wplywie perturbacji, tak aby dla malych krokow czasowych wszystkie
%czlony zmierzaly do zera}  NIE DZIALA

%\input impr-estm-srzed.tex  % nieudana proba lepszego oszacowania czlonow brezgowych uzywajac okresowosci


\begin{thebibliography}{ZCC}

%\bibitem[SEM]{SEM} Springer Encyplopedia of Mathematics, \begin{verbatim}http://www.encyclopediaofmath.org/index.php/Krylov-Bogolyubov_method_of_averaging %\end{verbatim}


\bibitem[AKN97]{AKN97} V. Arnold, V. Kozlov, and A. Neishtadt. \emph{Mathematical Aspects of Classical and Celestial Mechanics.}
Springer Verlag, Berlin, 1997.

\bibitem[BM]{BM} N.N. Bogolyubov and Y. A. Mitropol'skii, \emph{Asymptotic Methods in the Theory of Nonlinear Oscillations}, 2nd ed Gordon \&
Breach: New York, 1961

\bibitem[BZ]{BZ} N.N. Bogolyubov and D. N. Zubarev, \emph{An Asymptotic Approximation Method for a System with Rotating Phases and its Application to the Motion of a Charged Particle in a Magnetic Field}  Ukrain. Math. Zh. 7 (1955)

\bibitem[BL]{BL} A. Buica and  J. Llibre,  Averaging methods for finding periodic orbits via Brouwer degree, \emph{Bull. Sci. math.} 128 (2004) 7–22

\bibitem[CaZ]{CaZ} M.J. Capi\'nski and P. Zgliczy\'nski, Geometric proof for normally hyperbolic invariant manifolds, \emph{ J. Diff. Eq.}, 259(2015) 6215--6286

\bibitem[CMTZ]{CMTZ} J. Cyranka, P. B. Mucha, E. S. Titi, and P Zgliczy\'nski,
\emph{Stabilizing the Long-time Behavior of the Navier-Stokes Equations and Damped Euler Systems by Fast Oscillating Forces}, preprint
arXiv:1601.04612 (2016)

\bibitem[CyZ]{CyZ} J. Cyranka and P. Zgliczy\'nski, \emph{Stabilizing effect of large average initial velocity in forced dissipative PDEs invariant with respect to Galilean transformations}, preprint, arXiv:1407.1712 (2014)




%\bibitem[G]{G}
%S. Gerschgorin. \"{U}ber die Abgrenzung der Eigenwerte einer Matrix,
% \emph{ Bulletin de l'Acad\'{e}mie des Sciences de l'URSS. Classe des sciences math\'{e}matiques et na}, 6:749--754, 1931.

\bibitem[GH]{GH} J. Guckenheimer and P. Holmes. Nonlinear Oscillations, Dynamical Systems and Bifurcations of Vector Fields. Springer, New York, 1983.



\bibitem[HNW]{HNW} E.\ Hairer, S.P.\ N{\o}rsett and G.\ Wanner,
{\em Solving Ordinary Differential Equations I, Nonstiff  Problems}, Springer-Verlag, Berlin Heidelberg 1987.


\bibitem[Hale]{Hale} J.K. Hale,  Ordinary Differential Equations, Wiley, New York, 1969

\bibitem[He]{He} D. Henry, \emph{Geometric Theory of Semilinear Parabolic Equations}, in: Lecture Notes in Mathematics, vol. 840, Springer 1981



\bibitem[HPS]{HPS} M.W. Hirsch, C.C. Pugh and M. Shub, \emph{Invariant
manifolds}, Lecture Notes in Mathematics vol. 583, 1977


\bibitem[I]{I} A. Iserles, \emph{Three stories of high oscillation.}, Eur. Math. Soc.
Newsl. No. 87 (2013),

\bibitem[JS]{JS} A. Jorba, C. Simo, \emph{Effective Stability for Periodically Pertubed Hamiltonian Systems},  in 'Hamiltonian Mechanics
NATO ASI Series. Integrability and Chaotic Behavior' Volume 331, 1994, 245--252

\bibitem[KZ]{KZ}T. Kapela and P. Zgliczy\'nski, \textit{A Lohner-type algorithm for control systems and ordinary differential inclusions},
Discrete Cont. Dyn. Sys. B, vol. 11(2009), 365-385.

\bibitem[L]{L} J. Llibre, The Averaging Theory for Computing Periodic Orbits
Llibre, in
\textit{In Central Configurations, Periodic Orbits, and Hamiltonian Systems}, Birkhäuser Basel, CRM Barcelona, 2015.

\bibitem[MM]{MM} K. Mischaikow, M. Mrozek, \emph{Conley index.} Chapter 9 in Handbook of Dynamical Systems, vol 2, pp 393–460, Elsevier 2002

\bibitem[MS]{MS} M. Mrozek, R. Srzednicki, \emph{Topological Approach to Rigorous Numerics of Chaotic
Dynamical Systems with Strong Expansion of Error
Bounds}, Found. Comput. Math. 10 (2010)  191–220

\bibitem[Nei84]{Nei84} A. Neishtadt, \emph{The separation of motions in systems with rapidly rotating phase}, J.
Appl. Math. Mech. 48 (1984), 133-139.

\bibitem[OW]{OW} P. Oprocha, P. Wilczy\'nski, \emph{Distributional chaos via isolating segments}, DCDS B 8 (2007), 347--356


\bibitem[SVM]{SV} J. A. Sanders F. Verhulst and J. Murdock,
Averaging Methods in Nonlinear Dynamical Systems, Appl. Math. Sci., vol. 59, Springer, New York, 2007

\bibitem[Si]{Si} C. Simo, \emph{Averaging under Fast Quasiperiodic Forcing}, in 'Hamiltonian Mechanics
NATO ASI Series. Integrability and Chaotic Behavior' Volume 331, 1994,  13--34

\bibitem[S1]{S1} R. Srzednicki,
{\em Periodic and bounded solutions in blocks for time-periodic
nonautonomuous ordinary differential equations.\/}
Nonlin. Analysis, TMA., 1994, 22, 707--737

\bibitem[SW]{SW} R. Srzednicki,  K. W\'{o}jcik,
\emph{A geometric method for detecting chaotic dynamics\/},
J. Diff. Eq., 1997, 135, 66--82

\bibitem[Wa]{Wa} J. R. Ward Jr, \emph{Homotopy and Bounded Solutions of Ordinary Differential Equations}, J. Diff. Eq.,107 (1994), 428--445


\bibitem[WZ]{WZ} K. W\'ojcik, P. Zgliczy\'nski, \emph{Isolating segments, fixed point index and symbolic dynamics},
       J. Diff. Eq, 161, 245--288, (2000)

\bibitem[ZCC]{ZCC}P. Zgliczy\'{n}ski, \emph{Covering relations, cone
conditions and the stable manifold theorem}, J. Differential
Equations. 246(5), 1774-1819 (2009)



%\bibitem[ZGi]{ZGi}P. Zgliczy\'{n}ski and M. Gidea, \emph{Covering relations
%for multidimensional dynamical systems}, Journal of Differential
%Equations 202/1, 33-58 (2004)


\end{thebibliography}
\end{document}